\documentclass[12pt]{article}
\usepackage{latexsym,amsmath,amsfonts,amssymb}

\def \msk {\medskip}
\def \nin {\noindent}
\def \smin {\setminus}

\newtheorem{thm}{Theorem}

\newtheorem{prob}[thm]{Problem}

\newtheorem{quest}[thm]{Question}

\begin{document}

\title{The Cover Pebbling Number of Graphs}
\author{ Betsy Crull \thanks{
Department of Mathematics and Computer Science, Valparaiso University,
Valparaiso, Indiana 46383-6493. } \\
Tammy Cundiff $^{{\tiny *}}$\\
Paul Feltman $^{{\tiny *}}$\\
Glenn H. Hurlbert \thanks{
Department of Mathematics and Statistics, Arizona State University, Tempe,
AZ 85287-1804. email: hurlbert@asu.edu. Partially supported by National
Security Agency grant \#MDA9040210095. }\\
Lara Pudwell $^{{\tiny *}}$\\
Zsuzsanna Szaniszlo \thanks{
Department of Mathematics and Computer Science, Valparaiso University,
Valparaiso, Indiana 46383-6493. email: zsuzsanna.szaniszlo@valpo.edu. }\\
\ and \\
Zsolt Tuza \thanks{
Computer and Automation Institute of the Hungarian Academy of Sciences,
Budapest; and Department of Computer Science, University of Veszpr\'em,
Hungary. Partially supported by OTKA grant T-032969. }\\
}
\maketitle

\begin{abstract}
A pebbling move on a graph consists of taking two pebbles off of one vertex
and placing one pebble on an adjacent vertex. In the traditional pebbling
problem we try to reach a specified vertex of the graph by a sequence of
pebbling moves. In this paper we investigate the case when every vertex of
the graph must end up with at least one pebble after a series of pebbling
moves. The \textit{cover pebbling number} of a graph is the minimum number
of pebbles such that however the pebbles are initially placed on the
vertices of the graph we can eventually put a pebble on every vertex
simultaneously. We find the cover pebbling numbers of trees and some other
graphs. 
We also consider the more general problem where (possibly different) given
numbers of pebbles are required for the vertices.

\vspace{0.2 in}

\vspace{0.2 in}

\end{abstract}

\newpage

%
%

\newpage

%
%

\section{Introduction}

\label{Intro}

The game of pebbling was first suggested by Lagarias and Saks, and
introduced to the literature in a paper of Chung \cite{C}. A pebbling move
consists of taking two pebbles off of one vertex and placing one pebble on
an adjacent vertex. Given a graph $G$, a specified number of pebbles, and a
configuration of the pebbles on the vertices of $G$, the goal is to be able
to move at least one pebble to any specified target vertex using a sequence
of pebbling moves. The pebbling number $\pi (G)$ is the minimum number of
pebbles that are sufficient to reach any target vertex regardless of the
original configuration of the pebbles. In the present context it is
naturally assumed that \emph{all graphs considered are connected}. Moews
\cite{M} found the pebbling number of trees 
by using a clever path partition of the tree. For a survey of additional
results see \cite{H}.

In this paper we investigate the following question: How does the pebbling
problem change if instead of having a specified target vertex we need to
place a pebble simultaneously on every vertex of the graph? In some
scenarios this seems to be a more natural question, for example if
information needs to be transmitted to several locations of a network, or if
army troops need to be deployed simultaneously. We define the \textit{cover
pebbling number\/} ${\gamma }(G)$ to be the minimum number of pebbles needed
to place a pebble on every vertex of the graph using a sequence of pebbling
moves, regardless of the initial configuration. We establish the cover
pebbling number for several classes of graphs, including complete graphs,
paths, fuses (a fuse is a path with leaves attached at one end), and more
generally, trees. We also describe the structure of the largest
non-coverable configuration on a tree.

More generally, let a weight function $w$ be given that assigns an integer $%
w(v)$ to each vertex $v$ of $G$. We say that $w$ is \textit{positive} if $%
w(v)>0$ for all $v$. We define the \textit{weighted cover pebbling number\/}
${{\gamma}_w}(G)$ to be the minimum number $k$ ensuring that, from any
initial configuration with $k$ pebbles there is a sequence of pebbling moves
after which all the vertices $v$ simultaneously have $w(v)$ pebbles on them.
Our main result on trees in Section~\ref{Trees} determines ${{\gamma}_w}(T)$
for every tree $T$ and every positive weight function $w$.

\msk

Given a configuration $C$ of pebbles, we will use the following notation.
The \textit{size} $|C|$ of the configuration, denotes the number of pebbles
in $C$. The \textit{support} ${\sigma }(C)$ of the configuration is the set
of \emph{support vertices}, i.e. those on which there is at least one pebble
of $C $. The number of pebbles on $v$ in $C$ is denoted by $C(v)$ 
(hence, $v\in {\sigma }(C)$ if and only if $C(v)>0$). 
We call a configuration \textit{simple} if its support consists of a single
vertex. We say that a configuration is \textit{cover-solvable}, or simply
\textit{coverable} (resp. $w$-\textit{coverable}), if it is possible to
transport at least one pebble (resp. $w(v)$ pebbles) to every vertex $v$ of
the graph simultaneously (and \textit{non-coverable} otherwise). As is
customary, we denote the 
vertex set and edge set of $G$ by $V(G)$ and $E(G)$, respectively. 
If $G$ is of order $n$, we sometimes denote its vertices by $%
v_{1},v_{2},\dots ,v_{n}$.

%
%

\section{Preliminary Results}

\label{Results}

We begin with the cover pebbling number of the complete graph $K_{n}$ on $n$
vertices. Note that the pebbling number for $K_{n},$ $\pi (K_{n}),$ is $n$
(see \cite{H}).

\begin{thm}
\label{complete} ${\gamma}(K_n)=2n-1$.
\end{thm}

\noindent \textit{Proof.} If $2n-2$ pebbles are placed on vertex $v_n$, then
2 pebbles will be used to cover each of the $n-1$ other vertices. Thus no
pebbles will remain to cover $v_n$. Hence ${\gamma}(K_n)\ge 2n-1$.

Now suppose that at least $2n-1$ pebbles are placed on the vertices. We may
suppose that some vertex, say $v_n$, has no pebbles on it, otherwise the
graph is already covered. The pigeonhole principle says that some other
vertex has at least two pebbles on it; we use those to cover $v_n$. Since
there are now at least $2n-3$ pebbles among the remaining $n-1$ vertices,
induction says we can cover them (of course, ${\gamma}(K_1)=1$). Hence ${%
\gamma}(K_n)\le 2n-1$. \hfill ${\Box}$ \vspace{0.2 in}

A similar inductive proof works also for weighted covering, and yields the
following result. Denote the total weight by $|w|=\sum_v w(v)$ and define $%
\min w=\min_v w(v)$.

\begin{thm}
\label{complete-w} ${{\gamma}_w}(K_n)=2|w|-\min w$ for every positive weight
function\/ $w$.
\end{thm}

Next we find the cover pebbling number of the path $P_n$ on $n$ vertices $%
v_1,\ldots,v_n$, with $v_iv_{i+1}\in E$ for $1\le i<n$. Note that $%
\pi(P_n)=2^{n-1}$ (see \cite{H}).

\begin{thm}
\label{paths} ${\gamma}(P_n)=2^n-1$.
\end{thm}

\noindent \textit{Proof.} If $2^n-2$ pebbles are placed at vertex $v_n$,
then covering $v_1$ will use $2^{n-1}$ pebbles, covering $v_2$ will use $%
2^{n-2}$ pebbles, $\ldots$, and covering $v_{n-1}$ will use 2 pebbles. Then
no pebbles will remain to cover $v_n$. Hence ${\gamma}(P_n)\ge 2^n-1$.

Now suppose that at least $2^{n}-1$ pebbles are placed on the vertices. If
there are no pebbles on $v_{n}$ then we may use at most $2^{n-1}$ pebbles to
cover it, since $\pi (P_{n})=2^{n-1}$. By induction, the remaining $2^{n-1}-1
$ or more pebbles can cover $P_{n-1}$ (of course, ${\gamma }(P_{1})=1$). If
there are pebbles on $v_{n}$ then move as many of them as possible to $%
v_{n-1}$, leaving 1 or 2 on $v_{n}$. Either at least $2^{n-1}-1$
pebbles have been moved to $v_{n-1}$, or at most $2^{n-1}-2$ moves
have been made and at most two pebbles stay on $v_n$. In any case,
at
least $2^{n-1}-1$ pebbles remain on $P_{n-1}$. Again, induction shows that ${%
\gamma }(P_{n})\leq 2^{n}-1$. \hfill ${\Box }$ \vspace{0.2in}

Among all graphs on $n$ vertices, the complete graph has the smallest
pebbling number ($n$) and the path has the largest pebbling number ($2^{n-1}$%
). In both cases, we have ${\gamma }(G)=2\pi (G)-1$. While this might lead
one to guess that such a relation holds for all (connected) graphs, this
couldn't be farther from the truth. As the following theorem shows, the
ratio ${\gamma }(G)/\pi (G)$ is unbounded, even within the class of trees.
The subclass of \textit{fuses} is defined as follows. The vertices of $%
F_{l}(n)$ ($l\geq 2$ and $n\geq 3$) are $v_{1},\ldots ,v_{n}$, so that the
first $l$ vertices form a path from $v_{1}$ to $v_{l}$, and the remaining
vertices are independent and adjacent only to $v_{l}$. (The path is
sometimes called the \textit{wick}, while the remaining vertices are
sometimes called the \textit{sparks}.) For example, $F_{2}(n)$ is the star $%
S_{n}$ on $n$ vertices. The fact that ${\gamma }(S_{n})=4n-5$ serves as the
base case for the following result.

\begin{thm}
\label{fuses} ${\gamma}(F_l(n))=(n-l+1)2^l-1$.
\end{thm}

\noindent \textit{Proof.} Following the arguments for the path given above,
it is easy to see that so many pebbles are required of a simple
configuration sitting on $v_1$.

Likewise, induction on $n$ shows that so many pebbles suffice to cover the
fuse. Indeed, consider the cases whether or not $v_{1}$ has pebbles on it
and argue as was done for paths, above.

Regarding the base case $l=2$, we point out that $F_2(n)$ is the star on $n$
vertices, so we can let any leaf play the role of $v_1$. If all the pebbles
are on $v_2$ then we can cover the star easily. Otherwise, some leaf has at
least one pebble on it, and we label that vertex $v_1$. Now we pebble as
many as possible from $v_1$ to $v_2$, leaving 1 or 2 on $v_1$. Induction on
the number of leaves finishes the proof. \hfill ${\Box}$ \vspace{0.2 in}

We define the \textit{covering ratio} of $G$ to be ${\rho}(G)={\gamma}%
(G)/\pi(G)$. For a class $\mathcal{F}$ of graphs we define ${\rho}(\mathcal{F%
})=\sup_{G\in\mathcal{F}}{\rho}(G)$ if it exists, and ${\rho}(\mathcal{F}%
)=\infty$ otherwise. Thus, for the families $\mathcal{K}$ of complete graphs
and $\mathcal{P}$ of paths, we have ${\rho}(\mathcal{K})={\rho}(\mathcal{P}%
)=2$.

\begin{thm}
\label{ratio} Let $\mathcal{T}_n$ be the family of all trees on $n$
vertices. Then ${\rho}(\mathcal{T}_n)=\infty$.
\end{thm}

\noindent \textit{Proof.} Since $\pi(F_l(n))=2^l+n-l-1$ (see \cite{M}), we
see that, for $n=2^l+l$, ${\rho}(F_l(n)) > (n-l)2^l/(n-l+2^l) > (n-\lg n)/2$%
. \hfill ${\Box}$ \vspace{0.2 in}

%
%

\section{The Transition Digraph}

\label{Trans}

The main goal of this section is to prove that any sequence of pebbling
moves can be replaced by one which is cycle-free in a well-defined sense.
For this, we introduce the following concept.

\vspace{0.2 in}

\nin {\bf Definition.} Given a sequence $S$ of pebbling moves on graph $G$,
the \emph{transition digraph\/} is a \emph{directed multigraph\/} denoted $%
T(G,S)$ that has $V(G)$ as its vertex set, and each move $s\in S$ along edge
$v_iv_j$ (i.e., where two pebbles are removed from $v_i$ and one placed on $%
v_j$) is represented by one directed edge $v_iv_j$.

\msk

\begin{thm}
\label{cycle-free} Let\/ $S$ be a sequence of pebbling moves on\/ $G$,
resulting in a configuration\/ $C$. Then there exists a sequence\/ $S^*$ of
pebbling moves, terminating with a configuration\/ $C^*$, such that

\begin{enumerate}
\item On each vertex\/ $v$, the number of pebbles in\/ $C^*$ is at least as
large as that in\/ $C$, and

\item $T(G,S^*)$ does not contain any directed cycles.
\end{enumerate}
\end{thm}

\noindent \textit{Proof.} We apply induction on the number of directed
cycles in $T(G,S)$. The assertion is trivially true for every $S$ where this
number is zero. 

Let now $S$ be arbitrary, and consider the shortest prefix $S^{\prime}$ of $%
S $ that contains a directed cycle. That is, the last move in $S^{\prime}$
creates a cycle, say $C^{\prime}=v_1v_2\cdots v_k$, in $T(G,S^{\prime})$.
For $i=1,2,\dots,n$, let us denote by $d^-_i$ and $d^+_i$ the in-degree and
out-degree, respectively, of vertex $v_i$ in $T(G,S^{\prime})$. In the
initial configuration, each $v_i$ has to contain at least $2d^+_i - d^-_i$
pebbles, otherwise some move of $S^{\prime}$ could not be performed at $v_i$.

Let us consider the edge set $F^{\prime}= E(T(G,S^{\prime}))\smin %
E(C^{\prime})$. By the choice of $S^{\prime}$, this $F^{\prime}$ does not
contain any directed cycles, hence it contains a vertex $v_i$ of in-degree
zero. It means $d^-_i=0$ if $v_i\notin C^{\prime}$, and $d^-_i=1$ otherwise.
In the former case, $v_i$ initially has at least $2d^+_i$ pebbles and is
incident with precisely $d^+_i$ edges; while in the latter, the number of
pebbles at $v_i$ is at least $2d^+_i-1$ and that of its incident edges is
just $d^+_i-1$. In either case, $v_i$ has sufficiently many pebbles so that
the pebbling moves for all of its incident edges in $F^{\prime}$ are
feasible before any move belonging to $C^{\prime}$ has been performed. We
now rearrange $S^{\prime}$ to make all moves of $F^{\prime}$ involving $v_i$
at the beginning. Analogously, $F^{\prime}-v_i$ has a vertex $v_j$ of zero
in-degree in $F^{\prime}$, hence after the rearrangement of moves at $v_i$,
the moves of edges incident with $v_j$ are feasible completely before $%
C^{\prime}$. Eventually we obtain a rearrangement, say $S^{\prime\prime}$,
of $S^{\prime}$ where the moves of $C^{\prime}$ are performed at the very
end, and of course the concatenation of $S^{\prime\prime}$ and $S-S^{\prime}
$ terminates in configuration $C$. Now it is immediately seen that the
concatenation $S^+$ of $S^{\prime\prime}-C^{\prime}$ and $S-S^{\prime}$ is a
feasible sequence of moves that ends up with a configuration $C^+$ where the
vertices $v_1,\dots,v_k$ have one more pebble than in $C$, and the other
vertices have the same number of pebbles in $C$ and $C^+$. Since the number
of directed cycles in $T(G,S^+)$ is strictly smaller than that in $T(G,S)$,
the assertion follows by induction. \hfill ${\Box}$ \vspace{0.2 in}

%
%

\section{Trees}

\label{Trees}

In this section we determine the (weighted) cover pebbling number for an
arbitrary tree $T$. For $v\in V(T)$ define
\begin{equation*}
s(v)=\sum_{u\in V(T)}2^{d(u,v)}\ ,
\end{equation*}
where $d(u,v)$ denotes the distance from $u$ to $v$, and let
\begin{equation*}
s(T)=\max_{v\in V(T)}s(v)\ .
\end{equation*}
Analogously, if a positive weight function $w$ is given, we define
\begin{equation*}
s_w(v)=\sum_{u\in V(T)}w(u)\cdot 2^{d(u,v)}\
\end{equation*}
and
\begin{equation*}
s_w(T)=\max_{v\in V(T)}s_w(v)\ .
\end{equation*}
Clearly, for a \emph{simple\/} configuration sitting on $v$, $s_w(v)$
pebbles are necessary and sufficient to cover $T$. Thus ${{\gamma}_w}(T)\ge
s_w(T)$ for every $T$ and every positive $w$. We are going to prove that
this obvious lower bound is in fact tight.

\begin{thm}
\label{trees} For positive weight functions $w$ we have ${{\gamma}_w}%
(T)=s_w(T)$.
\end{thm}

\noindent \textit{Proof.} The theorem can be reformulated in the following
equivalent form:

\begin{quote}
\textsl{For every non-coverable configuration\/ $C$ there exists a
\underline{simple} non-coverable configuration\/ $C^*$ such that\/ $|C^*|=|C|
$.}
\end{quote}

The proof of this latter assertion is essentially induction, where we either
reduce the tree to another tree with fewer vertices or keep $T$ unchanged
but decrease the support ${\sigma}(C)$ of $C$ without making its size $|C|$
decrease.

We shall use the following terminology concerning a configuration $C$. We
say that a vertex $v$ is a

\begin{itemize}
\item D-vertex with demand $D(v)=w(v)-C(v)$ if $w(v)-C(v)>0$.

\item N-vertex (neutral) if $C(v)=w(v)$. Then we define $D(v)=0$.

\item S-vertex with supply $S(v)=C(v)-w(v)$ if $C(v)-w(v)>0$.
\end{itemize}

It is immediate by definition that every non-coverable configuration
contains at least one D-vertex.

\msk\nin{\bf Case 1.} {\sl $T=K_1$ or\/ $T=K_2$.}

These are trivial initial cases, handled already in the more general context
of Theorem~\ref{complete-w}.

\msk\nin{\bf Case 2.} {\sl Some leaf of\/ $T$ is not an S-vertex.}

Let $v$ be such a leaf, and let $u$ be its neighbor in $T$. We now delete $v$
from $T$ (with all its pebbles), and increase $w$ at $u$ to the value $%
w^{\prime }(u)=w(u)+2D(v)$. Keeping $w^{\prime }(x)=w(x)$ unchanged for all $%
x\notin \{u,v\}$, the configuration $C^{\prime }=C-v$ on the tree $T^{\prime
}=T-v$ with the weight function $w^{\prime }$ is coverable if and only if so
is $C$ on $T$ with $w$. This follows from Theorem~\ref{cycle-free}, which
implies that if $T$ is coverable, then there is a sequence of pebbling moves
where no pebble gets moved from $v$ to $u$. (To make $v$ properly covered,
we need to place at least $D(v)$ additional pebbles on it; and this requires
taking $2D(v)$ pebbles off of $u$.)

\msk\nin{\bf Case 3.} {\sl Every leaf of\/ $T$ is an S-vertex.}

For a given leaf $v=v_1$, define the path $v_1v_2\cdots v_m$ so that $v_m$
is the other leaf if $T$ is a path and is the only vertex of the path having
degree at least 3 in $T$ otherwise. In the latter case we call $v_m$ the
\textit{split} vertex of $v_1$. If there is a support vertex other than $v_1$
on this path, we call the one having minimum subscript the \textit{nearest
support} vertex of $v_1$.

Since $v_1$ is an S-vertex we can move $s_1=\lfloor S(v_1)/2\rfloor$ pebbles
to $v_2$. Moreover, if $s_1>w(v_2)-C(v_2)$ then we can further transmit $%
s_2=\lfloor (s_1+C(v_2)-w(v_2))/2\rfloor$ pebbles to $v_3$, and so on. For a
vertex $v_k$ on this path we say that $v_1$ \textit{supplies} $v_k$ if at
least one of the pebbles from $v_1$ can reach $v_k$ in this way. There are
three possibilities for $v_1$, namely, $v_1$ supplies its split vertex, $v_1$
supplies its nearest support vertex, or $v_1$ supplies neither of these. We
consider these possibilities in reverse order.

\msk%
\nin{\bf Subcase A.} {\sl Some leaf supplies neither its split nor its nearest
support vertices.}

We follow a similar argument as in Case~2. Let $v_{1}$ be such a leaf and
let $k$ be the largest subscript so that $v_{1}$ supplies $v_{k}$ (then $k<m$
and $v_{i}$ is not a support vertex for any $2\leq i\leq k$). Let $C^{\prime
}$ and $w^{\prime }$ be the restrictions of $C$ and $w$ to $T^{\prime
}=T-\{v_{1},\ldots ,v_{k}\}$, except that $w^{\prime
}(v_{k+1})=w(v_{k+1})+2D^{\prime }$, where $D^{\prime }=w(v_{k})-s_{k-1}$ is
the resulting demand on $v_{k}$ after being supplied by $v_{1}$. Then $%
C^{\prime }$ is non-$w^{\prime }$-coverable on $T^{\prime }$, and since $%
|T^{\prime }|<|T|$ there is a simple non-$w^{\prime }$-coverable
configuration of size $|C^{\prime }|$ on $T^{\prime }$. This yields a non-$w$%
-coverable configuration $C^{\prime \prime }$ of size $|C|$ on $T$ that sits
on two vertices. If $T$ has at least three leaves then some leaf is not an
S-vertex and we are done by Case~2. Otherwise $T$ is a path and ${\sigma }%
(C^{\prime \prime })=\{v_{1},v_{n}\}$. Non-$w$-coverability now means that $%
v_{n}$ can supply $v_{k}$ with strictly fewer than $D^{\prime }$ pebbles.
Finally we test if $k-1\geq n-k$. If so, then for every $j$ in the range $%
k\leq j\leq n$, $d(v_{1},v_{j})\geq d(v_{j},v_{k})$. Thus, defining $C^{\ast
}(v_{n})=0$ and $C^{\ast }(v_{1})=C^{\prime }(v_{1})+C^{\prime }(v_{n})=|C|$%
, we obtain a simple non-coverable configuration, as required. If $k-1<n-k$
we do the opposite.

\msk\nin{\bf Subcase B.} {\sl Some leaf supplies its nearest support vertex.}

Let $v_1$ be such a vertex and $v_k$ its nearest support vertex (then $%
v_i\notin{\sigma}(C)$ for $1<i<k$). We define $C^{\prime}(v_k)=0$ and $%
C^{\prime}(v_1)=C(v_1)+C(v_k)$, keeping $C^{\prime}$ identical to $C$ on
every other vertex. Then $|C^{\prime}|=|C|$, $|{\sigma}(C^{\prime})|<|{\sigma%
}(C)|$, and $C^{\prime}$ is non-coverable whenever $C$ is, because the
supply from $v_1$ yields fewer pebbles on $v_k$ in $C^{\prime}$ than in $C$.

\msk\nin{\bf Subcase C.} {\sl Every leaf supplies its split vertex.}

By Subcase~B we may assume that no leaf supplies its nearest support vertex.
There must be some vertex $v$ that is the split vertex for two different
leaves (indeed, choose any leaf and let $v$ be any vertex of degree at least
3 at farthest distance from it -- the two leaves past $v$ witness this).
Label these leaves $v_{1}$ and $v_{\ell }$ so that $P=v_{1}\cdots
v_{m}\cdots v_{\ell }$ is the unique path between them and $v=v_{m}$. Recall
that $v_{i}$ is not a support vertex for any $1<i<\ell $ and that both $%
v_{1} $ and $v_{\ell }$ supply $v_{m}$. Let us denote by $s_{m}$ their total
supply for $v_{m}$.

If $s_{m}>w(v_{m})$, then $P$ can supply $T-P$ with $s^{\prime }=\lfloor
\frac{1}{2}(s_{m}-w(v_{m}))\rfloor $ pebbles (at most); and otherwise it
needs to receive at least $s^{\prime \prime }=w(v_{m})-s_{m}$ pebbles from $%
T-P$. In both cases we consider the problem restricted to $P$, where $%
w(v_{i})$ is kept unchanged for all $i\neq m$, and $w(v_{m})$ is modified to
$s_{m}+1$. This configuration on $P$ is non-coverable. Thus, according to
Subcase~A, the $C(v_{1})+C(v_{\ell })$ pebbles can be placed on one vertex ($%
v_{1}$ or $v_{\ell }$), keeping $P$ non-coverable. It follows that the
modified configuration, too, either supplies $T-P$ with at most $s^{\prime }$
pebbles or needs to receive at least $s^{\prime \prime }$ pebbles from $T-P$%
. In either case, the new configuration on $T$ is non-coverable and has at
least one D-vertex leaf, thus we are done by Case~2.

\hfill ${\Box}$ \vspace{0.2 in}

From this proof we see that a non-coverable configuration of maximum size
can be assumed to be simple. The next result shows that the single support
vertex must be an end of a longest path. (This is the case even for weight
functions $w$ where the longest paths are not of maximum weight.)

\begin{thm}
\label{config} Given a tree $T$ and a positive weight function $w$, let $C$
be a non-coverable simple configuration of maximum size, with ${\sigma}%
(C)=\{v\}$. Then $v$ is a leaf of a longest path in $T$.
\end{thm}

\noindent \textit{Proof.} Since ${{\gamma}_w}(T)=s_w(v)$ for some $v$, we
need to show that the maximum value of $s_w(v)$ is attained only on some
endpoints of the longest path(s) of $T$. We are going to prove something
stronger: \emph{every\/} longest path has at least one endpoint $x$ whose $%
s_w(x)$ is larger than $s_w(u)$ for \emph{every\/} $u$ which is not the
endpoint of some longest path.

Suppose first that $T$ is just a path $v_1v_2\cdots v_n$. Consider any
internal vertex $v_k$ ($1<k<n$). We compare the partial sums $s^-=\sum_{1\le
i<k} w(v_i)\cdot 2^{d(v_i,v_k)}$ and $s^+=\sum_{k< i\le n} w(v_i)\cdot
2^{d(v_i,v_k)}$. If $s^-\le s^+$, then $s_w(v_{k-1})>s_w(v_k)$; and if $%
s^-\ge s^+$, then $s_w(v_{k+1})>s_w(v_k)$. Thus, $s_w(k)$ can never be
largest.

Suppose next that $T$ is a tree with precisely three leaves. Applying the
previous idea, from any non-leaf vertex we can move to one of its neighbors
and find there a larger value of $s_w$. Hence, let $v,v^{\prime},v^{\prime%
\prime}$ be the three leaves, and suppose that the longest path $P$ in $T$
is the one connecting $v^{\prime}$ with $v^{\prime\prime}$. We need to show $%
s_w(v)<\max\,\{s_w(v^{\prime}),s_w(v^{\prime\prime})\}$. Let $u$ be the
unique degree-3 vertex of $T$. We have $d(u,v)<d(u,v^{\prime})$ and $%
d(u,v)<d(u,v^{\prime\prime})$ (for otherwise the $v$--$v^{\prime}$ path or
the $v$--$v^{\prime\prime}$ path were at least as long as the $v^{\prime}$--$%
v^{\prime\prime}$ path, contrary to the assumption on $v$). From this it is
easily seen that for every vertex $x$, at least one of $d(v^{\prime},x)$ and
$d(v^{\prime\prime},x)$ is at least $d(v,x)+1$. Consequently, $%
s_w(v^{\prime})+s_w(v^{\prime\prime})>2s_w(v)$, i.e.\ $s_w(v)$ cannot be
largest.

Finally, let $T$ be a tree with more than three leaves. Let $P$ be one of
its longest paths, $v^*$ a leaf that does \emph{not\/} belong to any longest
path of $T$, and $v\neq v^*$ a leaf not on $P$ (but maybe on some other
longest path of $T$). We apply the transformation on $v$ as described in
Case~2 of the proof of Theorem~\ref{trees}. This modification keeps the
function $s_w$ unchanged on all vertices of $T-v$, moreover $P$ remains a
longest path and $v^*$ does not become the endpoint of any longest path in $%
T-v$. Thus, by induction on the number of vertices, $s_w$ is larger on some
endpoint of $P$ than on $v^*$. This completes the proof. \hfill ${\Box}$
\vspace{0.2 in}

%
%

\section{Open Problems}

\label{Open}

There are several natural problems and questions to ask.

\begin{prob}
Find ${\gamma }(G)$ for other graphs $G$, for example cubes, complete $r$%
-partite graphs, etc.
\end{prob}

\begin{quest}
Is it true for all graphs $G$ that at least one of the largest
non-coverable configurations on $G$ is simple?
\end{quest}

\begin{prob}
Find classes of graphs $\mathcal{F}$ whose covering ratio ${\rho}(\mathcal{F}%
)$ is bounded.
\end{prob}

\begin{quest}
Can the question, ``Is ${\rho}(G)\le k$?'' be answered efficiently?
\end{quest}

These questions extend to positive weight functions in a natural way. Let us
note, however, that the situation drastically changes when ``positive'' is
replaced by ``nonnegative'' for $w$. This fact is already shown by the
complete graph $K_n$ ($n\ge 3$) where only one vertex is required to be
covered, which corresponds to the weights $1,0,0,\dots,0$. Here the unique
maximal non-coverable configuration has the pebble distribution $%
0,1,1,\dots,1$, in striking contrast to the case where $w>0$ and all pebbles
may be concentrated on a suitably chosen single vertex. Such considerations
must be tackled in order to pursue the \textit{weighted pebbling number} of
a graph $G$, defined as $\pi_{\mathrm{w}}(G) = \max_w {\gamma}_w(G)$, where
the maximum is taken over all nonnegative weight functions $w$ of size $|w|=%
\mathrm{w}$. The pebbling number $\pi(G)$ is the case $\mathrm{w}=1$.

\begin{prob}
Find $\pi_{\mathrm{w}}(T)$ for any tree $T$ and weight \textrm{w}.
\end{prob}


%
%

%
%
%

%
%

\end{document}